\catcode`\^^@9
\catcode`\9
\catcode`\∙10
\documentclass[12pt]{article}
\usepackage[T2A]{fontenc}
\usepackage[english]{babel}
\usepackage{upref,xspace}
\usepackage{amssymb,amscd}
\usepackage[dvips]{graphicx}
\usepackage{enumerate}
\usepackage{amsmath,amsthm}
\setkeys{Gin}{draft=false}
\def\bey#1\eey{\begin{equation}#1\end{equation}}
\def\ben#1\een{\begin{equation*}#1\end{equation*}}
\def\bst#1\est{\begin{split}#1\end{split}}

\newcommand{\rom}[1]{\textrm{\mdseries\upshape#1}}

\let\cal\mathcal
\let\Bbb\mathbb

\let\ge \geqslant

\let\le \leqslant
\let\leq\leqslant
\def\qh{\enskip}

\def\qq{\qquad}
\def\qqq{\qquad\quad}
\def\qu{\quad}
\def\tss#1{\ \text{#1}\ }
\def\fract#1#2{#1/#2}
\let\SUBSECTION\subsection
\renewcommand{\subsection}[1]
{\refstepcounter{subsection}\SUBSECTION*{\thesubsection. #1}}
\let\SUBSUBSECTION\subsubsection
\renewcommand{\subsubsection}[1]
{\refstepcounter{subsubsection}\relax
\SUBSUBSECTION*{\bfseries\thesubsubsection. #1}\relax
}
\newcommand{\mes}{\mbox{mes\,}}
\newcommand{\supp}{\mbox{supp\,}}
\newcommand{\Z}{{\Bbb Z}}
\newcommand{\R}{{\Bbb R}}
\newcommand{\N}{{\Bbb N}}
\newcommand{\Cvar}{{\Bbb C}}
\newcommand{\ve}{\varepsilon}
\newcommand{\ity}{\infty}
\newcommand{\ti}{\tilde}
\newcommand{\be}{\beta}
\newcommand{\de}{\delta}
\newcommand{\De}{\Delta}
\newcommand{\la}{\lambda}
\newcommand{\al}{\alpha}
\newcommand{\vfi}{\varphi}
\newcommand{\si}{\sigma}

\newcommand{\Bl}{\Bigl(}
\newcommand{\Br}{\Bigr)}
\newcommand{\Bm}{\Bigr|}
\newcommand{\eps}{\varepsilon}
\newcommand{\inz}{\int\limits_{2^n\le |z|\le 2^{n+1}}}
\newcommand{\ilm}{\int\limits}
\renewcommand{\Im}{\mathop{\rm Im}}
\renewcommand{\Re}{\mathop{\rm Re}}
\theoremstyle{plain}%обеспечивает в автомате курсив
\newtheorem{theorem}{Theorem}
\newtheorem{corollary}{Corollary}
\newtheorem{lemma}{Lemma}

\newtheorem{pheo}{Theorem}
\renewcommand{\thepheo}{\Alph{pheo}}
\newtoks\thehthmvar
\newtheorem*{thmvar}{\the\thehthmvar}
\newenvironment{thm}[1]
 {\thehthmvar{#1}\begin{thmvar}}
 {\end{thmvar}}
\theoremstyle{definition}%обеспечивает
\newtoks\thehRemark
\newtheorem*{Remarkvar}{\the\thehRemark}
\newenvironment{Remark}[1]
 {\thehRemark{#1}\begin{Remarkvar}}
 {\end{Remarkvar}}
\numberwithin{equation}{section}
\theoremstyle{plain}
\begin{document}
\begin{center} \large\bf APPROXIMATION
OF SUBHARMONIC FUNCTIONS\end{center}
\begin{center}{\large I. Chyzhykov} \end{center}
%\address{
\begin{center}
\textit{Faculty of Mechanics and Mathematics, \\
Lviv National University, Lviv, Ukraine%}
\\
ichyzh@lviv.farlep.net }\end{center}
%\thanks{1991 {\em Mathematics Subject Classification}.
%              Primary 30A05; Secondary 30D20, 30E10.}
%}
%\author{Игорь Чижиков$^{0}$
%\thanks{Partly supported by  project .}
%}
%\keywords{  }
%\thanks{The research is partially supported
%by the scholarship of the Queen Jadwiga foundation, Jagellonian
%University (Krak\'ow, Poland)}

%\runningheads{I.~Chyzhykov} {Approximation of subharmonic
%functions}

\begin{abstract}
In certain classes of subharmonic functions $u$ on $\mathbb{C}$ distinguished in terms of lower bounds for the
Riesz measure of $u$,  a sharp estimate is obtained for the 
rate of  approximation  by 
functions of the form $\log |f(z)|$, where $f$ is an entire
function. The results complement and generalize those recently
obtained by Yu. Lyubarskii and Eu. Malinnikova.
\end{abstract}
\date{}
%\maketitle
%\thanks Исследование частично поддержано
%стипендией фонда королевы Ядвиги  Ягеллонского университета
%(Краков, Польша)\endthanks
\section{Introduction}
We use the standard notions of the subharmonic function theory
\cite{Hay}. Let us introduce some notation. Let $D_z(t)=\{\zeta\in
\Cvar: |\zeta-z|<t\}$,   $z\in \Cvar$, $t>0$. For a subharmonic
function $u$ in $\Cvar$ we write $B(r,u)=\max \{ u(z): |z|=r\}$,
$r>0$, and define the order $\rho[u]$ by the equality
$\rho[u]=\limsup\limits_{r\to+\ity} \log  B(r,u)/\log  r$. Let
also $\mu_u$ denote the Riesz measure associated with the
subharmonic function $u$,    $m$ be the plane Lebesgue measure.
The symbol $C$ with indices stands for some positive constants. If
$u=\log  |f|$, where $f$ is an entire function with the zeros
$\{a_k\}$, then $\mu_u=\sum\limits_k n_k \de(z-a_k)$, where $n_k$
is the multiplicity of the zero  $a_k$, $\de(z-a_k)$ is the Dirac
function concentrated at the point~$a_k$. However, the class of
functions subharmonic in $\Cvar$ is broader than those of the form
$\log  |f|$, where $f$ is an entire function. Since it is often
easier to construct a subharmonic function rather than an entire
one with desired asymptotic properties, a natural problem arises
of approximation of subharmonic functions by the logarithms of the
moduli of entire functions. Apparently, V.\ S.\ Azarin \cite{Az}
was the first to investigate
 this problem in the
general form in the class of subharmonic in the plane functions
and of finite order of growth. The results cited below have
numerous applications in the function theory and the potential
theory (see, e.g., %\cite{LS}--\cite{Ch1}).
%^^@было:\cite{LS--Ch1}).
%надо:
\rom{[3--6]}).%номера взяты из распечатки

In 1985  R.\ S.~Yulmukhametov \cite{Yul} obtained the following
remarkable result. {\it For each function $u$ subharmonic in
$\Cvar$ and of order $\rho\in (0,+\ity)$, and $\alpha>\rho$ there
exist an entire function $f$ and a set  $E_\al\subset \Cvar$ such
that
\begin{equation}
\big|u(z)-\log   |f(z)|\big|\le C_\al \log  |z|, \quad  z\to\ity,
\qh z\not\in E_\al, \label{e1.1}
\end{equation}
and $E_\al$ can be covered by a family of disks $D_{z_j}(t_j)$,
$j\in \N$, satisfying  the estimate  $\sum_{|z_j|>R}
t_j=O(R^{\rho-\al})$, $(R\to+\ity)$.}
In the special  case when  a subharmonic function $u$ is
homogeneous $\log |z|$ in estimate (\ref{e1.1})  can be replaced
by $O(1)$
 \cite{LS,Yul2}.
An integral metric allows us to drop an exceptional set when we
approximate subharmonic functions of finite order. Let  $\|\cdot
\|_q$ be the norm in the space $L^q(0, 2\pi)$,
%^^@неправильно:
%^^@$$Q(r,u)=\cases O(\log   r), \; r\to+\infty, &\rho[u]<+\infty,\\
%O(\log   r +\log   n(r,u)), \; r\to+\infty,&\\
%\quad  r\not\in E, \mes E<+\infty, &
%\rho[u]=+\infty .\endcases$$
%правильно (надо \begin-, \end{cases}):
$$ Q(r,u)=\begin{cases} O(\log  r), \; r\to+\infty,
&\rho[u]<+\infty,\\ {\aligned[b] O(\log  r +\log  n(r,u)), \;
r\to+\infty,&\\ \quad  r\not\in E, \mes E<+\infty,&
\endaligned}
&
\rho[u]=+\infty.
\end{cases}
$$ A.\ A.~Gol'dberg and M.\ O.\ Hirnyk  have proved \cite{GG} that
for an arbitrary  subharmonic  function
  $u$ there exists an entire function  $f$ such that
$$ \big\| u(re^{i\theta}) -\log   |f(re^{i\theta})|\big\|_q
=Q(r,u), \quad r>0, \qh q>0. $$ From the recent result of Yu.\
Lyubarskii and Eu.\ Malinnikova \cite{LM} it follows that
integrating of an approximation rate by the plane measure allows
us to drop the assumption that a subharmonic function $u$ is of
finite order of growth and to obtain sharp estimates.

%\begin{pheo} \cite{LM} %{ Theorem А.\/}
%{\it  Пусть $u(z)$ --- субгармоническая в $\Cvar$ функция.
% Тогда  для любого $q>1/2$ существуют число
%  $R_0>0$ и  целая функция $f$ такие, что
%\begin{equation}
%\frac 1{\pi R^2} \int\limits_{|z|<R}
%         |u(z)-\log  |f(z)||dm(z) <q \log   R, \quad R> R_0.
%\label{eq:in1e}
%\end{equation}
%\label{th:in2}
%}\end{pheo}
\begin{thm}{\refstepcounter{pheo} Theorem  \thepheo\ \cite{LM}} %{ Theorem А.\/}
\label{th:in2} Let  $u(z)$ be a subharmonic function in $\Cvar$.
Then, for each $q>1/2$, there exist
  $R_0>0$ and an entire function  $f$ such
that
\begin{equation}
\frac 1{\pi R^2} \int\limits_{|z|<R}
         \big|u(z)-\log |f(z)|\big|dm(z) <q \log  R, \quad R> R_0.
\label{eq:in1e}
\end{equation}
\end{thm}
An example constructed in \cite{LM} shows that we cannot  take
$q<1/2$ in estimate (\ref{eq:in1e}).
In connection with Theorem A the following M.\ Sodin's question,
which is a more precise form of Question 1 from
\cite[p.~315]{Sod}, is  known:

\medskip
\noindent {\bf Question.} {\it Given a subharmonic function $u$ in
$\Bbb C$, do there exist  an entire function $f$ and a constant
 $\al\in [0,1)$ such that
\begin{equation}
 \int\limits_{|z|<R}
         \big|u(z)-\log |f(z)|-\al \log |z|\big|\,dm(z) =O(R^2), \quad R> R_0.
\label{e:sh}
\end{equation}}
^^@\begin{Remark}{Remark 1} Question 1 from \cite[p.~315]{Sod}
corresponds to the case $\al=0$. The mentioned example from
\cite{LM} implies the negative answer to this question.
\end{Remark}
On the other hand, the restrictions onto the Riesz measure from
below allow us to sharpen estimate  (\ref{eq:in1e}).
\begin{thm}{\refstepcounter{pheo}Theorem \thepheo\ \cite{LM}}%{\bf Теорема Б
Let $u(z)$ be a subharmonic in $\Cvar$ function. If, for some
$R_0>0$ and $q>1$
\begin{equation}
\mu_u(\{z:R<|z|\leq qR\})>1, \quad R>R_0,
\label{e1.2}
\end{equation}
then there exists an entire function $f$ satisfying
 $$
\displaystyle \sup\limits_{R>0} R^{-2} \int\limits_{|z|<R}
         \big|u(z)-\log  |f(z)|\big|\,dm(z) <\infty.
$$
%\label{eq:in1c}
In addition, for every $\eps >0$ there exists a set $ E_\eps
\subset \Cvar$ such that $$ \displaystyle\limsup_{R\to \infty}
          {m(\{z\in E: |z|<R\})} R^{-2} <\eps$$
%\label{eq:in1d}
and
\begin{equation}
u(z)-\log  |f(z)| = O(1), \qu  z\not\in E_\eps, \qh z\to \infty.
\label{e:in1}
\end{equation}
\end{thm}
As we can notice, there is a gap between the statements of
Theorems A and B. The following question arises: how can one
improve the estimate of the left hand side of (\ref{eq:in1e}),
whenever (\ref{e1.2}) fails to hold? The answer to this question
is given by Theorem \ref{t:1}.

Let $\Phi$  be the class of slowly varying functions $\psi\colon
[1,+\ity)\to (1,+\ity)$ (in particular, $\psi(2r)\sim \psi(r)$ for
$r\to+\ity$).
%such that $r/\psi(r) \nearrow +\infty$ as $r\to+\ity$
\begin{theorem}\label{t:1}
Let $u$ be a subharmonic in $\Cvar$ function, $\mu=\mu_u$. If for
some $\psi\in \Phi$ there exists a constant $R_1$ satisfying the
condition
\begin{equation}
(\forall R>R_1): \mu(\{ z:R<|z|\le R\psi(R)\})>1,
\label{e1.3}
\end{equation}
then there exists an entire function $f$ such that $(R\ge R_1)$
\begin{gather}
\int\limits_{|z|<R}
         \big|u(z)-\log  |f(z)|\big|\,dm(z) = O(R^2 \log   \psi (R)).
\label{e1.4}
\end{gather}
\end{theorem}

\begin{corollary}
In the conditions of Theorem \ref{t:1}, for arbitrary $\ve>0$
there exist $K(\ve)>0$ and a set
 $E_\ve\subset \Cvar$ such that
\begin{equation}
\limsup_{r\to\infty} \frac{m(E\cap \overline{D}_r)}{r^2} <\ve
\label{e1.8}
\end{equation}
and \begin{equation}\label{e1.9} \big|u(z)-\log  |f(z)|\big|\le
K_\ve \log  \psi(|z|), \quad z\not\in E_\ve.
\end{equation}
\end{corollary}
The following example demonstrates that estimate (\ref{e1.4}) is
exact in the class of subharmonic functions satisfying
(\ref{e1.3}). This is indicated by Theorem~\ref{t:2}.

For $\vfi\in \Phi$, let
\begin{equation} \label {ex2.1}
u(z)= u_\vfi(z)=\frac12 \sum_{k=1}^{+\ity} \log   \Bigl| 1-\frac
z{r_k}\Bigr|,
\end{equation}
where $r_0=2$, $r_{k+1}=r_k \vfi(r_k)$, $k\in \N\cup \{ 0\}$.
Thus, $\mu_{u}$ satisfies condition (\ref{e1.3}) with
$\psi(x)=\vfi^2(x)$.
\begin{theorem} \label{t:2}
%\it
Let  $\psi\in\Phi$ be such that $\psi(r)\to +\infty$
$(r\to+\infty)$. There exists no entire function $f$ for which
\begin{equation*}
\int\limits_{|z|<R}
         \big|u_\psi(z)-\log  |f(z)|\big|\,dm(z) = o(R^2 \log   \psi (R)), \quad
R\to\infty.
\end{equation*}
\end{theorem}
It follows from the proved Theorem $2'$ that the answer to the
formulated above M.~Sodin's question is negative. It was M. Hirnyk
who pointed out the question as well as the fact that the negative
answer follows from the example constructed above.

%\begin{theorem} \label{t:2'} %получается теорема 3!
\medskip \noindent
{\bf Theorem $\bf 2'$.} {\it Let $\sigma$ be an arbitrary positive
continuous, defined on $[1,+\infty)$ function and $ \si(t)\to 0$
$(t\to\infty)$, $\psi(R)=\exp\Bigl\{ \int_1^R \frac{\si(t)}t \,
dt\Bigr\}$. There is no entire function $f$ and constant
$\al\in [0,1)$ for which
\begin{equation*}
\int\limits_{|z|<R} \big|u_\psi(z)-\log  |f(z)| -\al\log
|z|\big|\,dm(z) = o\Bigl(R^2 \log \psi(R)\Bigr), \quad R\to\infty.
\end{equation*}
} ^^@\begin{Remark}{Remark 2} The growth of the expression
$\int_1^R \frac{\si(t)}t \, dt$ for $R\to +\infty$ is restricted
only by the condition $\int_1^R \frac{\si(t)}t \, dt=o(\log  R)$.
\end{Remark}
^^@\begin{Remark}{Remark 3} The author does not know whether it is
possible to improve estimate (\ref{e1.8}) of the exceptional set
for (\ref{e1.9}). In \cite{Ch2} there are obtained the sharp
estimates of the exceptional set outside of which relation
(\ref{e1.9}) holds, for some class of subharmonic functions
satisfying some additional restriction onto the Riesz measure.
\end{Remark}
\section{Proof of Theorem \ref{t:1}}
\SUBSECTION*{2.1. Partition of  measures} It is appeared that
in order to have a ``good'' approximation of a
 subharmonic function by the logarithm of the  modulus of
an entire function we need a ``good'' approximation of the
corresponding Riesz measure by a discrete one. The only
restriction on Riesz' measure,  defined on the Borel sets in the
plane, is finiteness  on the compact sets. The following theorem
on a partition of measures is the principal step in proofs of the
theorems cited above.
\begin{pheo} \label{t:ge}
 \it Let $\mu$ be a measure in  $\R^2$ with compact support,
$\rom\supp\mu \subset \Pi$, and $\mu (\Pi)\in \N$, where $\Pi$ is
a square. Suppose, in addition, that for any line $L$ parallel to
either side of the square  $\Pi$, there is at most one point $p\in
L$ such that
\begin{equation} \label{e:2a}
 0<\mu(\{p\}) <1
%\quad \textrm{ и }  \quad
\tss{and} \mu(L\setminus \{p\})=0 ,
\end{equation}
Then there exist a system of rectangles $\Pi_k\subset \Pi$ with
sides parallel to the sides of $\Pi$, and measures $\mu_k$ with
the following properties\rom:
\begin{itemize}
\item[1)]  $\rom\supp\mu_k \subset \Pi_k$\rom;
\item[2)] $\mu_k(\Pi_k)=1$, $\sum_k\mu_k =\mu$\rom;
\item[3)]  the interiors of the  convex
hulls of the supports of  $\mu_k$ are pairwise disjoint\rom;
\item[4)] the ratio
of length of the sides for rectangles $\Pi_k$ lies  in the
interval $[1/3,3]$\rom;
\item[5)] each point of the plane belongs
to the interiors of at most 4  rectangles~$\Pi_k$.
\end{itemize}
\end{pheo}
Theorem \ref{t:ge} was proved  by R.\ S.\ Yulmukhametov
\cite[Theorem 1]{Yul} for absolutely continuous measures (i.e.
$\nu$ such that $m(E)=0 \Rightarrow \nu(E)=0$). In this case
condition (\ref{e:2a}) holds automatically. As it is shown by
D.~Drasin~\cite[Theorem 2.1]{Dr}, Yulmukhametov's proof works if
we replace the condition of continuity by condition
 (\ref{e:2a}). We can drop
condition (\ref{e:2a}) rotating the initial square~\cite{Dr}.
Though  it is noted in \cite{LM} that Theorem  \ref{t:ge} is valid
even without hypothesis (\ref{e:2a}), the author does not know a
proof of this fact. Moreover, proving Theorem 2.1 \cite{Dr} (a
variant of Theorem \ref{t:ge}) condition (\ref{e:2a}) is used
essentially. In this connection one should mention   a note by
A.F.~Grishin and S.V.~Makarenko \cite{GM}, where a two dimensional
variant of Theorem \ref{t:ge} under the condition that the measure
does not load lines  parallel to the coordinate axes.
^^@\begin{Remark}{Remark 4} In the proof of Theorem \ref{t:ge}
\cite{Dr} rectangles $\Pi_k$  are obtained  by partition of  given
rectangles, starting from $\Pi$,  into smaller rectangles in the
following way. Length of a smaller side of an initial  rectangle
coincides with length of a side of obtained rectangle, and length
of the other side of obtained rectangle is not less than a third
part and not greater than  two third parts of length of the other
side for the initial rectangle. Thus the following form of
Theorem~\ref{t:ge}, which will be used, holds.
\end{Remark}
\begin{theorem} \label{t:5}
Let $\mu$ be  a measure in $\R^2$ with compact support, $\rom\supp
\mu \subset \Pi$, $\mu (\Pi)\in 2\N$, where $\Pi$ is a rectangle
with the ratio $b_0/a_0=l_0\in [1, +\infty)$ of length
 $a_0, b_0$, $(a_0\le b_0)$ of the sides.
If, in addition, condition \rom(\ref{e:2a}\rom) holds, then there
exist a system of rectangles  $\Pi_k\subset \Pi$ with sides
parallel to the sides of $\Pi$ and measures $\mu_k$ with the
following properties\rom:
\begin{itemize}
\item[1)] $\rom\supp \mu_k \subset \Pi_k$\rom;
\item[2)] $\mu_k(\Pi_k)=2$, $\sum_k\mu_k =\mu;$
\item[3)]  the interiors of the  convex hulls of the
supports of $\mu_k$ are pairwise disjoint\rom;
\item[4)] the ratio $b_k/a_k$ of  length $a_k, b_k$ $(a_k\le  b_k)$
of the sides for the rectangle $\Pi_k$ lies in the interval
 $[1,l_0]$, moreover, if
$l_k>3$, then $a_k=a_0$\rom;
\item[5)] each
point of the plane belongs to the interiors of at most 4
rectangles~$\Pi_k$.
\end{itemize}
\end{theorem}
As it is noted in \cite{LS}, the idea of  partition into
rectangles with mass 2 is due to A.~F.~Grishin.
In order to apply Theorem \ref{t:5} we need the following lemma
(see also Lemma 2.4~\cite{Dr}.)
\begin{lemma}
Let $\nu$ be  a locally finite measure in  $\Cvar$. Then in any
neighborhood  of the origin there exists a point
 $z'$ with the following properties\rom:
\begin{itemize}
\item[a)]  on each line $L_\alpha$ going through   $z'$
there is at most one point  $\zeta_\alpha$ such that
 $\nu(\{\zeta_\alpha\})>0$, moreover  $\nu(L_\al\setminus \{\zeta_\al\})=0;$
\item[b)]  on each circle $C_\rho$ with center
 $z'$ there exists at most one point $\zeta_\rho$ such
that
 $\nu(\{\zeta_\rho\})>0$, moreover $\nu(C_\rho\setminus
\{\zeta_\rho\})=0$.
\end{itemize}
\end{lemma}
We  give a simple example for illustration.
Let $\nu(z)=\sum_{n\in \N} \delta(z-n)$. Then $\nu(\Bbb
R)=+\infty$. We can take  any point of the disk $\{ z: |z|<1\}$
with an nonzero imaginary part as $z'$.

\begin{proof}[^^@\bf\upshape Proof of Lemma  1]
Let $B_n=\{z: 2^n< |z|\le 2^{n+1}\}$, $n\in \N$, $B_0=\{z : |z|\le
2\}$ and  $\nu_n=\nu\bigr|_{B_n}$. Since $\nu$ is a locally finite
measure, $\nu_n(\Cvar)=\nu_n(B_n)<+\infty$. There  is an at most
countable set $\zeta_{nk}$ of points such that
$\nu_n(\{\zeta_{nk}\})>0$.  Therefore the set $E_1=\bigcup _n
\bigcup_k \{\zeta_{nk}\}$ is at most countable. Given a pair of
points from $E_1$ we consider the straight line going through
these points, and the  middle perpendicular to the segment
connecting these points. All these lines cover some set
 $A\subset \Cvar$, $m(A)=0$.
Let $z_1\in \Cvar\setminus A$. By our construction, an arbitrary
straight line going through the point $z_1$ contains at most  one
point with positive mass. The same is true for an arbitrary circle
with center $z_1$. We define
$\nu_n'=\nu_n-\sum_{\nu_n(\{\zeta\})>0}
\nu_n(\{\zeta\})\delta_\zeta$ $(n\in \Z_+)$,
$\delta_\zeta(z)=\delta(z-\zeta)$. Then, for any $z\in \Cvar$, we
have $\nu'_n (\{z\})=0$ $(n\in \Z_+)$. Since the intersection of
two  different  circles (straight lines)  is  empty or a point or
two points, and $\nu_n(\Cvar)<+\infty$, by  countable additivity
of $\nu_n$, there exists  at most countable  set of circles and
straight lines  with  positive  $\nu_n$-measure. The union   $F_n$
of all these straight lines and  centers of  circles has zero
area. If now $z'\in \Cvar \setminus (A\cup \bigcup_{n\in \Z_+}
F_n)$, then for any  $n\in \Z_+$ the measure $\nu_n$ of  any
circle with center $z'$ as well as that of a straight line going
through $z'$ is equal to zero. Hence, their  $\nu$-measure  equals
 zero. Finally, by  the countable additivity  of the plane
measure
 $m(A\cup \bigcup_{n\in \Z_+}
F_n)=m(A)+\sum_n m(F_n)=0$. Thus, any point $z'\in(\Cvar \setminus
(A\cup \bigcup_{n\in \Z_+} F_n))\cap U $ where $U$ is a given
neighborhood of the origin has required properties.
  % \qed
\end{proof}
Taking into account the proved lemma one may assume that
properties a) and b) of Lemma 1 hold for the origin. We follow the
scheme of the proof from \cite{LM}, assuming that
$\psi(x)\nearrow+\infty$ $(x\to+\infty)$, because otherwise
Theorem 1 is equivalent to Theorem B. Without loss of generality,
one may assume that $u(z)$ is harmonic in a neighborhood of $z=0$.
Otherwise, choose arbitrary $a>0$ such that $n(a)\le N\le n(a+0)$
for some $N\in \N$ and define the measure $\nu$ to be equal $\mu$
in $D_0(a)$ and to contain the part $\mu\bigr|_{\{z:|z|=a\}}$ so
that $\nu(\overline{D_0(a)})=N$. Then $\mu-\nu \equiv 0$ in
$D_0(a)$ and, instead of $u$, we consider the function  $\tilde
u(z)=u(z)- \int_{\Cvar} \log  |z-\zeta|\, d\nu(\zeta)$,  and the
quantity $\int_\Cvar \log  |z-\zeta|\, d\nu(\zeta)- N\log  |z|$ is
bounded whenever $z\to +\infty$. Therefore, without loss of
generality we assume that $R_0=\sup\{ r>0: \supp \mu \cap
D_0(r)=\varnothing\}>1$. Define $\Psi_1(R)=R\psi(R)$,
$\Psi_n(R)=\Psi_1(\Psi_{n-1}(R))$ for $n\in \N$, $\Psi_0(R)\equiv
R$, $R>1$. We define by induction measures $\mu_k^{(j)}$, $j\in
\{1,2,3\}$, and a sequence $(R_k)$, $k\in \N$. Suppose that
 $\mu_l^{(j)}$ is already defined for $l<k$, and $R_l$ for $l\le k$.
Let $$ Q_k=\{ \zeta\in \Cvar: R_k\le |\zeta|\le \Psi_1(R_k)\},
\quad \mu_k^-=\Bigl(\mu- \sum _{j=1}^{k-1}
(\mu_j^{(1)}+\mu_j^{(2)}+\mu_j^{(3)})\Bigr)\biggr|_{Q_k}. $$
If $\mu_k^{-} (Q_k)<2$, define $\mu_k^{(1)}(Q_k)\equiv 0$,
$\mu_k^{(2)}=\mu_k^-$, $\mu_k^{(3)} \equiv 0$.
If $\mu_k^{-} (Q_k)\ge 2$, represent  $\mu_k^- $ as the sum
$\mu_k^{(1)}+\mu_k^{(2)}$, so that $\mu_k^{(1)}(Q_k)=
2[\mu_k^-(Q_k)/2]$, where $[a]$ denotes the integer part of $a$.
If  $\mu_k^{(2)}(Q_k)\ge 1$, we define $\mu_k^{(3)}(Q_k)=0$ and
$R_{k+1}=\Psi_1(R_k)$. Otherwise, let $$ R_{k+1}= \inf\bigl\{ R>
R_k: \mu(\{\Psi_1(R_k)<|z|\le R\})\ge 1\bigr\}, $$ and by
$\mu_k^{(3)}$ we mean the sum of the restriction of $\mu$ onto
$\{\zeta: \Psi_1(R_k)<|\zeta|< R_{k+1}\}$ and $\tilde \mu$, where
$\tilde \mu$ is the part of the restriction
$\mu\bigr|_{\{\zeta:|\zeta|=R_{k+1}\}}$ such that
$\mu_k^{(3)}(\Cvar)=1$.
By the construction, we have for all $k\in \N$:
\begin{itemize}
\item[1)] $\supp \mu_k^{(1)} \subset Q_k$, $\mu_k^{(1)} (Q_k)\in 2\Z_+$;
\item[2)] $\Psi_1(R_k)\le R_{k+1} \le \Psi_2(R_k)$;
\item[3)] $\supp (\mu_k^{(2)}+\mu_k^{(3)}) \subset \{ \zeta: R_k\le |\zeta |
\le R_{k+1} \}$;
\item[4)] $1\le (\mu_k^{(2)}+\mu_k^{(3)})\bigl(\{ \zeta: R_k\le |\zeta |
\le R_{k+1} \}\bigr) \le 2$.
\end{itemize}
Let $\mu^{(1)}=\sum_{j=1}^{+\infty} \mu_j^{(1)}$,
$\mu^{(2)}=\sum_{j=1}^{+\infty} (\mu_j^{(2)}+\mu_j^{(3)})$. From
properties 3) and 4) it follows that
$\mu^{(2)}(\overline{D_0(R)})=O(\log  R)$ $(R\to +\infty)$,
therefore $u_2(z)=\int_{\Cvar} \log  \bigl| 1-\frac z\zeta \bigr|
\, d\mu^{(2)}(\zeta)$ is a subharmonic function in $\Cvar$. Let
$u_1(z)=u(z)-u_2(z)$. Then $\mu_{u_1}=\mu^{(1)}$. We will
approximate $u_1$ and $u_2$ separately. It suffices to prove that
\begin{equation} \label{e2.2}
I_n\stackrel{\rm def}=\int\limits_{2^n\le |z|\le 2^{n+1}}
\big|u(z)-\log  |f(z)|\big| \, dm(z) = O( 4^n \log  \psi (2^n)),
\quad n\to+\infty.
\end{equation}
Indeed, let $R\in [2^n, 2^{n+1})$. Then from (\ref{e2.2}) it
follows that
%\begin{gather*}
\begin{align*}
&\frac{\int\limits_{|z|<R} |u(z)-\log  |f(z)||\, dm(z)}{R^2} \le
\frac{\sum\limits_{k=0}^n \int\limits_{2^k \le |z|\le 2^{k+1}}
|u(z)-\log  |f(z)||\, dm(z)+ O(1) }{4^{n}}
\\
&\qq\le\frac{ C_1\sum\limits_{k=0}^  n 4^k \log  \psi(2^k) }{4^n}
\le 4C_1  \log  \psi (2^n).
%\end{gather*}
\end{align*}
\SUBSECTION*{2.2. Approximation of $u_2(z)$} By the construction,
$\mu^{(2)}(\Cvar)=+\infty$. Define
 $T_n=\sup\{R>0: \mu^{(2)}(\overline{D_0(R)})\le 5n \}$, $A_n=\{ \zeta: T_n \le |\zeta|\le
T_{n+1}\}$. Let $(A_n, \mu_n)$ be a partition of the measure
$\mu^{(2)}$ such that $\mu^{(2)}=\sum_{k=1}^{+\infty} \mu_k$,
$\supp \mu_n\subset A_n$, $\mu_n(A_n)=5$. Define $r_n$ by the
equalities $\log  r_n=\frac 15 \int_{A_n} \log  |\zeta| \,
d\mu_n(\zeta)$, $n\in \N$, and consider the formal product
 $$ f_2(z)=\prod_{n=1}^{+\infty} \Bigl( 1- \frac
z{r_n}\Bigr)^5. $$ From property 4) of the measures $\mu_k^{(j)}$
it follows that
\begin{equation} \label{e2.1}
\Psi_1(T_n)\le T_{n+1} \le \Psi_6(T_n).
\end{equation}
Since $r_{n+1}/r_{n-1}\ge T_{n+1}/T_{n}\ge \psi(T_n) \to+ \infty$
$(n\to+\infty)$, the function $f_2$ is entire.
Let
%\begin{equation}
\begin{align}
d_k(z)&\equiv \int\limits _{A_k} \Bl \log   \Bigl| 1-\frac z \zeta
\Bigr| -\log   \Bigl| 1- \frac z{r_k} \Bigr| \Br \,
d\mu_{k}(\zeta)\nonumber
\\
&=\int\limits _{A_k} \Bl \log \Bigl| 1-\frac \zeta z \Bigr| -\log
\Bigl| 1- \frac {r_k}z \Bigr| \Br \, d\mu_{k}(\zeta). \label
{e2.2'}
%\end{equation}
\end{align}
Here we make use the choice of $r_k$. Fix $n\in \N$, and let $2^n
\in [T_N, T_{N+1})$. Then for $k\ge N+2$ and $\zeta \in A_k$  we
have $|\zeta|\ge T_k \ge |z| {T_k}/{T_{N+1}}$, $ r_k \ge
|z|{T_k}/{T_{N+1}}$, consequently $ \bigl| \log | 1-\frac z \zeta
| \bigr|\le 2 |z|/|\zeta|$, $\bigl|\log  | 1-\frac z
{r_k}|\bigr|\le 2 |z|/r_k$. Therefore,
%\begin{gather}
\begin{align}
&\inz \sum_{k=N+2}^{\infty}  |d_k(z)|\, dm(z)\le  \nonumber
\inz \sum_{k=N+2}^{\infty}  20 \frac{T_{N+1}}{T_k} dm(z)
\\
&\qq\le C_2 \frac {T_{N+1}}{T_{N+2}} 4^n = o(4^n), \quad n\to \infty.
\label{e2.3}
%\end{gather}
\end{align}
Similarly, $|d_k(z)|\le 20 T_{k+1}2^{-n} $ for $k\le N-2$. Using
(\ref{e2.2'}), we obtain
\begin{equation} \label{e2.4}
\inz \sum_{k=1}^{N-2}  |d_k(z)|\, dm(z)\le
o(4^n), \quad n\to \infty.
\end{equation}
Estimate $\int_{2^n\le |z|\le  2^{n+1}} |d_k(z)|\, dm(z)$ for
$k\in \{N-1, N, N+1\}$. By the definition,
%$$
\begin{align*}
&\int\limits_{2^n\le |z|\le  2^{n+1}} |d_{N+1}(z)|\, dm(z)
\\
&\qq=\int\limits_{2^n\le |z|\le  2^{n+1}} \biggl|
\int\limits_{A_{N+1}} \Bigl( \log  \Bigl| 1-\frac z\zeta \Bigr|
-\log  \Bigl| 1-\frac z{r_{N+1}} \Bigr| \Bigr) \,
d\mu_{N+1}(\zeta)\biggr| dm(z).
\end{align*}
%$$
\newcommand{\izn}{\int\limits_{2^n\le |z|\le  2^{n+1}}}
\renewcommand{\Bl}{{\Bigl(}}
\renewcommand{\Br}{{\Bigr)}}
\renewcommand{\Bm}{{\Bigl|}}
If $|\zeta|\ge 2|z|$, $\zeta\in A_{N+1}$, we have $\bigl| \log  |
1-\frac z\zeta |\bigr| \le \log  2$. Otherwise, $T_{N+1}\le
|\zeta| < 2|z|< 2^{n+2} \le 4T_{N+1}$. Therefore, applying
Fubini's theorem and changing the variables, $T_{N+1}\eta=\zeta$,
$T_{N+1}\xi=z$, we obtain
%\begin{gather}\nonumber
\begin{align}
&\int\limits_{2^n\le |z|\le  2^{n+1}}  \int\limits_{A_{N+1}}
%\biggl|\log  \Bigl| 1-\frac z\zeta \Bigr| \biggr|
\Bm \log  \Bm 1-\frac z\zeta \Bm\Bm \, d\mu_{N+1}(\zeta) dm(z)
\nonumber
\\
&\le  2\izn dm(z)+\izn %\nonumber
%\\
%&\qq\qq\times
\int\limits_{T_{N+1}\le |\zeta| \le 4T_{N+1}} \Bm \log  \Bm
1-\frac z\zeta \Bm\Bm \, d\mu_{N+1}(\zeta)dm(z) \nonumber
\\
&\qu=O(4^n) \nonumber
\\
&\qq+ T_{N+1}^2 \nonumber
\\
&\qq\qu\times\int\limits_{1\le |\eta| \le 4}  \;
\int\limits_{{2^n}/{T_{N+1}}\le |\xi|\le {2^{n+1}}/{T_{N+1}}} \Bm
\log  \Bm 1-\frac \xi\eta \Bm\Bm \, dm(\xi)
\,d\mu_{N+1}(T_{N+1}\eta)\nonumber
\\
&\qu=O(4^n)+ O(4^n) \int\limits_{1\le |\eta| \le 4} \;
\int\limits_{|\xi|\le 2} \Bm \log  \Bm 1-\frac \xi\eta \Bm\Bm \,
dm(\xi) \,d\mu_{N+1}(T_{N+1}\eta). \label{e2.5}
%\end{gather}
\end{align}
However, elementary calculations demonstrate that for $1\le \eta
\le 4$ we have $$ \int\limits_{|\xi|\le 2} \Bm \log  \Bm 1-\frac
\xi\eta \Bm\Bm \, dm(\xi)\le C_3. $$ Taking into account that
$\mu_{N+1}(\Cvar)=2$, from (\ref{e2.5}) we deduce
\begin{equation}\label{e2.6}
\izn \int\limits_{A_{N+1}} \Bm \log  \Bm 1-\frac z\zeta \Bm\Bm
d\mu_{N+1}(\zeta)\, dm(z)=O(4^n), \quad n\to+\infty.
\end{equation}
Similarly, for $r_{N+1} \ge 2^{n+2} \ge 2|z|$ we have $\big| \log
|1- z/{r_{N+1}}|\big|\le \log  2$, otherwise, $T_{N+1}\le r_{N+1}
< 2^{n+2}$, and consequently  $(T_{N+1}\xi =z)$,
%\begin{gather}\nonumber
\begin{align}\nonumber
&\izn \Bm \log  \Bm 1-\frac z{r_{N+1}} \Bm\Bm\, dm(z)
\\
&\qq=  O(4^n)+ O\biggl(4^n \int\limits_{1/2\le |\xi|\le 1} \Bm
\log  \Bm 1-\frac \xi\eta \Bm\Bm\, dm(\xi) \biggr)= O(4^n), \quad
n\to+\infty. \label{e2.7}
%\end{gather}
\end{align}
Now, we are going to estimate the integral of $|d_N(z)|$ over the
annulus $\{2^n\le|z|\le 2^{n+1}\}$. For $\vfi(x)\equiv x/\psi(x)$
we have
%\begin{gather*}
\begin{align*}
&\int\limits_{A_N} \izn \Bm \log  \Bm \frac {z-\zeta }\zeta
\Bm\Bm\, dm(z)\, d\mu_N(\zeta)
\\
&\qq\le \int\limits_{2^{n-1}\le |\zeta|\le 2^{n+2}} \biggl(
\int\limits_{D_\zeta(\vfi(|\zeta|))}+\int\limits_{\begin{substack}
{ 2^n\le |z|\le 2^{n+1}\\ z\not\in D_\zeta(\vfi(|\zeta|))
}\end{substack}} \biggr) \Bm \log  \Bm \frac {z-\zeta }\zeta
\Bm\Bm\, dm(z)\, d\mu_N(\zeta)
\\
&\qq\qu+\biggl( \int\limits_{T_N\le |\zeta| \le 2^{n-1}}
+\int\limits_{2^{n+2}\le|\zeta|\le T_{N+1}} \biggr)
\int\limits_{2^n\le |z|\le 2^{n+1}} \Bm \log  \Bm \frac {z-\zeta
}\zeta \Bm\Bm\, dm(z)\, d\mu_N(\zeta)
\\
&\qq\qu\equiv I_{N,1}+I_{N,2}+I_{N,3}+I_{N,4}.
%\end{gather*}
\end{align*}
Integrating by parts and taking into account the relationship
$\vfi(x)=o(x)$ $(x \to +\infty)$, we obtain
%\begin{gather*}
\begin{align*}
I_{N,1} &=\int\limits_{2^{n-1} \le|\zeta| \le 2^{n+2}} \biggl(
2\pi \int\limits_0^{\vfi(|\zeta|)} \log  \frac{|\zeta|}s s\,
ds\biggr)\, d\mu_N(\zeta)
\\
&=2\pi \int\limits_{2^{n-1} \le |\zeta|\le 2^{n+2}} \biggl(
\frac{\vfi^2(|\zeta|)}2 \log  \psi(|\zeta|)+
\int\limits_0^{\vfi(|\zeta|)} s\,ds\biggr) d\mu_N(\zeta)
\\
&=(\pi +o(1))\int\limits_{2^{n-1}\le|\zeta| \le 2^{n+2}}
\vfi^2(|\zeta|)\log \psi(|\zeta|)\, d\mu_N(\zeta)= o(4^n), \quad
n\to+\infty;
%\end{gather*}
\end{align*}
%\begin{gather*}
\begin{align*}
I_{N,2} &\le \int\limits_{2^{n-1}\le |\zeta|\le 2^{n+2}}
\int\limits_ {\begin{substack} { 2^n\le |z|\le 2^{n+1}\\ z\not\in
D_\zeta(\vfi(|\zeta|)) }\end{substack}} \log
\frac{|\zeta|}{\vfi(|\zeta|)}\, dm(z)\, d\mu_n(\zeta)
\\
&\le \log  \psi(2^{n+1}) 4^{n+1} \pi \mu_N(\{ 2^{n-1}\le
|\zeta|\le 2^{n+2}\}) = O(4^n \log  \psi(2^n)),
\\
&\phantom{\le \log  \psi(2^{n+1}) 4^{n+1} \pi \mu_N(\{ 2^{n-1}\le
|\zeta|\le 2^{n+2}\}) = O(4 \log } n\to+\infty.
%\end{gather*}
\end{align*}
In the expressions for $I_{N,3}$ and $I_{N,4}$ the relationships
between $\zeta$ and $z$ yield respectively
%\begin{gather*}
\begin{align*}
\Bm \log  \Bm \frac{z-\zeta}\zeta \Bm\Bm &\le  \log  \frac{
2^{n+1} +T_N}{T_N} \le \log  \Bl 1+ 2\frac {T_{N+1}}{T_N}\Br
\\
&\le \log  \psi(T_N) +O(1),
\\
\Bm \log  \Bm \frac{z-\zeta}\zeta \Bm\Bm &\le 2\Bm \frac z\zeta
\Bm \le 1.
%\end{gather*}
\end{align*}
Therefore, $I_{N,3} +I_{N,4}= O(4^n \log  \psi(2^n) )$
$(n\to+\infty)$. In a similar way one can obtain the estimate
 $$ \izn \Bm \log  \Bm
1-\frac z{r_N} \Bm \Bm \, dm(z)=O(4^n \log  \psi(2^n)), \quad
n\to+\infty. $$ Combining the estimates for $I_{N,j}$, $j\in\{1,
\dots, 4\}$, and the latter inequality, we obtain
\begin{equation} \label{e2.8}
\izn |d_N(z)|\, dm(z)= O\big(4^n \log  \psi (2^n) \big), \quad
n\to+\infty.
\end{equation}
From definition (\ref{e2.2'}) it follows that
%$$
\begin{align*}
&\izn |d_{N-1}(z)|\, dm(z)
\\
&\qq \leqslant \izn  \int\limits_{A_{N-1}} \Bm \log  \Bm 1-\frac
\zeta z\Bm - \log  \Bm 1- \frac {r_{N-1}}z\Bm \Bm \,
d\mu_{N-1}(\zeta)\, dm(z).
\end{align*}
%$$
Analogously as the integral of $d_{N+1}(z)$ was estimated, for
$|\zeta|\le |z|/2$, $\zeta\in A_{N-1}$ we have $|\log  |
1-z/\zeta||\le \log  2$, therefore
%\begin{gather*}
\begin{align*}
&\izn \int\limits_{A_{N-1}} \Bm \log  \Bm 1-\frac \zeta z \Bm \Bm
d\mu_{N-1}(\zeta)\, dm(z)
\\
&\qu\le \izn \int\limits _{2^{n-1} \le |\zeta|\le T_N} \Bm \log
\Bm 1-\frac \zeta z \Bm \Bm d\mu_{N-1}(\zeta)\, dm(z)+ O(4^n)
\\
&\qu\le \int\limits _{T_N/2\le |\zeta|\le T_N} \izn \Bm \log  \Bm
1-\frac \zeta z \Bm \Bm dm(z)\, d\mu_{N-1}(\zeta) + O(4^n)
\\
&\qu\le T^2_N \ilm_{\frac12 \le |\eta| \le 1}  \ilm _{2^n/T_N\le
|\xi|\le 2^{n+1}/T_N} \Bm \log  \Bm 1-\frac \eta \xi \Bm \Bm
dm(T_{N}\xi)\, d\mu_{N-1}(T_N\eta)  +O(4^n)
\\
&\qq=O(T^2_{N}) +O(4^n) =O(4^n),
\quad n\to +\infty.
%\end{gather*}
\end{align*}
Finally, from (\ref{e2.6})--(\ref{e2.8}) and the latter
relationship we obtain
\begin{equation} \label{e2.9}
\izn \big|u_2(z)-\log  |f_2(z)|\big| \, dm(z) =O(4^n \log
\psi(2^n)), \quad n\to+\infty.
\end{equation}

\SUBSECTION*{2.3. Approximation of the function $u_1(z)$.} Let us
recall that $Q_k=\{\zeta: R_k\le |\zeta|\le R_k \psi(R_k)\}$,
$\supp \mu_k^{(1)} \subset Q_k$, $\mu_k^{(1)} (Q_k)\in 2\Z_+$. Let
 $$
P_k=\log  Q_k=\{s=\si+it: \log  R_k \le \si \le  \log  R_k +\log
\psi(R_k), 0\le t< 2\pi\}. $$ Denote by $l_k$ the ratio of the
larger side of the rectangle $P_k$ to the smaller one. For  $k\ge
k_0$ we have $l_k={\log  \psi(R_k)}/{(2\pi)}>1$. We consider the
measure $\mu_k^{(1)}$. If $\mu_k^{(1)}(\{p\})\ge 2$ at some point
$p$, we subtract from this measure the measure $\tilde\mu_k^{(1)}$
which is equal to $2[ \mu_k^{(1)}(\{p\})/2]$ at every such  $p$.
The measure $\tilde \mu^{(1)}=\sum_k \tilde \mu_k^{(1)}$ is
discrete, integer-valued, and finite on the compacta in $\Cvar$.
By the Weierstrass theorem, there exists an entire function $f_3$
with zeros of the corresponding multiplicity on the support of the
measure $\tilde \mu^{(1)}$, which is the  union of at most
countable set of isolated points. We have $\mu_{\log |f_3|}=\tilde
\mu^{(1)}$. For every  $k\in \N$ the measure $\mu^k
=\mu_k^{(1)}-\tilde \mu^{(1)}\bigr|_{Q_k}$ satisfies the condition
$\mu^k (\{p\})<2$  at every point $p\in Q_k$. According to the
choice of the origin, on the rays emanating from it as well as on
the circles centered at the origin, there is at most one point $p$
such that $0< \mu^k (\{p\})<2$, and at the same time the
$\mu^k$-measure of the remaining part of either ray or circle
equals zero. Under these conditions, the measures $\nu^k$ defined
by the conditions $d\nu^k(s)=d\mu^k(e^s)$, $s\in P_k$, (i.~e.
$\nu^k(S)=\mu^k(\exp S)$ for every Borel set $S$), satisfy the
conditions of Theorem \ref{t:5} with $\Pi=P_k$, $k\in \N$.
Applying this theorem, we obtain a system of rectangles $P_{km}$
and measures $\nu_{km}$, $1\le m \le N_k$, respectively. We have
$\nu_{km}(P_{km})=2$, every point $s\in \{ s: 0\le |\Im s|\le
2\pi\}$ belongs to the interiors of at most four rectangles
$P_{km}$. Enumerate $(P_{km}, \nu_{km})$ by the natural numbers,
in arbitrary way. As the result, we obtain $(P^{(k)}, \nu^{(k)})$
with $\nu^{(k)}(P^{(k)})=2$, $\supp \nu^{(k)} \subset P^{(k)}$;
then also the rest of properties from the formulation of Theorem
\ref{t:5} holds. ^^@\newcommand{\om}[1]{\omega_l^{(#1)}}
\newcommand{\ze}[1]{\zeta_l^{(#1)}}
\newcommand{\pl}{P^{(l)}}
\newcommand{\nl}{\nu^{(l)}(\omega)}
\renewcommand{\L}{{\cal L}}
\newcommand{\w}{{\omega}}
\newcommand{\diam}{\mathop{\rm diam}}
^^@Let $\om1$, $\om2$ be the solutions of the system of equations
\begin{equation}
\label{e3.1}
%^^@\cases \displaystyle \om1+\om2 =\int\limits_{P^{(l)}} \w d\nl, \\
%\displaystyle (\om1)^2+(\om2)^2 =\int\limits _{\pl} \w^2 d\nl,
%\endcases
\left\{{\aligned
\om1+\om2 &=\int\limits_{P^{(l)}} \w d\nl, \\
(\om1)^2+(\om2)^2 &=\int\limits _{\pl} \w^2 d\nl,
\endaligned}\right.
\end{equation}
and
\begin{equation}\label{e3.2}
\w_l=\frac12 \int\limits_{\pl} \w d\nl
\end{equation}
the center of mass of $\pl$, $l\in \N$. From (\ref{e3.1}) we find
$\om1+\om2=2\w_l$, $$ \int\limits_{\pl} \w^2d\nl =(\om1)^2 +
(2\w_l -\om1)^2= 2(\om1 -\w_l)^2+2\w_l^2. $$ Therefore, taking
into account (\ref{e3.2}), we obtain $(j\in \{1,2\})$
%$$
\begin{align*}
&|\om{j}-\w_l|=\Bigl| \frac12 \int\limits_{\pl} \w^2 d\nl - \w_l^2\Bigr|
^{\frac12}
\\
&\qq=\Bm \frac 12 \int\limits_{\pl} (\w -\w_l)^2 d\nl \Bm ^{\frac12 }\le
\diam \pl\equiv d_l.
\end{align*}
%$$
Since $\w_l\in \pl$, we obtain
\begin{equation} \label{e3.3}
\sup_{\w\in \pl} |\w-\om{j}|\le 2d_l, \qh j\in\{1,2\},  \quad
\sup_{\w\in \pl} |\w-\w_{l}|\le d_l.
\end{equation}
Let $\ze{j}=\exp\{\om{j}\}$,
%\begin{gather} \nonumber
\begin{align} \nonumber
V(z)&=\sum_l \int\limits_{Q^{(l)}} \Bl \log  \Bm 1-\frac z\zeta
\Bm -\frac 12 \log  \Bm 1-\frac z{\ze1} \Bm -\frac 12 \log  \Bm
1-\frac z {\ze2} \Bm\Br \, d\mu^{(l)}(\zeta) \nonumber
\\
&=\sum_l \int\limits_{\pl} \Bl \log  |1-ze^{-\omega}| -\frac 12
\log  \Bm 1-ze^{-\om1}\Bm -\frac12 \log  \Bm 1-ze^{-\om2}\Bm\Br\,
d\nl \nonumber
\\
&\equiv \sum_l \Delta _l(z).
\label{e3.4}
%\end{gather}
\end{align}
In the assumption that series (\ref{e3.4}) absolutely converges,
we have to prove that
 $$ \inz |V(z)| \, dm(z)=O(4^n \log
\psi(2^n)). $$
%^^@(23.03.04 20:22)
{\tolerance4000\looseness1 \noindent Let $\L^+$ be the set of $l$s
for which $Q^{(l)}\subset D_0(2^{n-2})$, and $\L^-$ for which
$Q^{(l)} \subset \{ z: |z|>2^{n+1}\}$, $\L^0=\N \setminus
 (\L^- \cup \L^+)$.
Denote $L(\w)=\log  (1-ze^{-\w})$. For $l\in \L^-\cup \L^+$,
$\w\in \pl$, $2^n\le |z| \le 2^{n+1}$, the function $L(\w)$ is
analytic. We will need the following equalities
%\begin{gather}\nonumber
\begin{align}\nonumber
&L(\w)-L(\om1)=\int\limits_{\om1}^\w L'(s)\, ds
\\
&\qq= L'(\w_l)(\w-\om1) +
\int\limits_{\om1}^\w L''(s)(\w-s)\, ds \nonumber
\\
&\qq=L'(\w_l)(\w-\om1)
\nonumber
\\
&\qqq
+\frac12 L''(\w_l)(\w-\om1)^2
+\frac 12\int\limits_{\om1}^\w L'''(s)(\w-s)^2\, ds.
\label{e3.5}
%\end{gather}
\end{align}\par}
^^@\noindent Using the first equality (\ref{e3.5}), we obtain
%\begin{gather} \nonumber
\begin{align} \nonumber
|\Delta_l(z)|&=\Bm \Re \int\limits_{\pl} \Bl L(\w)-L(\om1)-\frac12
(L(\om2)-L(\om2)\Br \, d\nl\Bm \nonumber
\\
&\le \int\limits_{\pl} \biggl| \int\limits_{\om1}^\w \frac {z}{e^s-z} ds\biggr|
\, d\nl+\frac 12 \int\limits_{\pl} \biggl| \int\limits_
{\om1}^{\om2} \frac {z}{e^s-z} ds\biggr| \, d\nl \nonumber
\\
&\le |z| \sup_{s\in \pl}\frac 1{|e^s-z|}
\biggl(\int\limits_{\pl} (|\w-\om1|+\frac12 |\om1-\om2|)
\, d\nl \biggr) \nonumber
\\
&\le \frac {3d_l|z|}{\inf_{s\in\pl} |e^s-z|}.
\label{e3.6}
%\end{gather}
\end{align}
Let $l\in \L^-$. The rectangle $\pl$ is of the form $\{s=\si+it:
\si_l^-\le \si\le \si_l^+, 0\le t_l^- \le t\le t_l^+ \le 2\pi\}$,
i.e. the lengths of its sides are equal to $\si_l^+-\si_l^-$ and
$t_l^+-t_l^-$. Let $\la_l$ be the ratio of the maximal of these
numbers to the minimal one $(\la_l\ge 1)$. According to condition
4) of the statement of Theorem \ref{t:5}, for $\la_l>3$ we have
$t_l^+-t_l^-=2\pi$, $\si_l^+-\si_l^-=2\pi\la_l$. First, we
estimate $|\Delta_l(z)|$ for which $\la_l>3$. Then $d_l=2\pi
\sqrt{1+\la_l^2}\le 2\log  \psi(e^{\si_l^-})$, because $\log
R_k\le \si_l^-< \si_l^+\le \log  R_k+\log  \psi(R_k)$ for some
$k\in \N$. From (\ref{e3.6}) it follows that
 $$ \inf\limits_{s\in \pl}
|e^s-z| \ge \inf_{s\in \pl} \frac{e^{\Re s}}2 \ge
\frac{e^{\si_l^-}}2, \quad   l\in \L^-. $$ But $$
\int\limits_{{\si_l^-}}^{{\si_l^+}} e^{-s} \, ds =e^{-\si_l^-}(1 -
e^{-\si_l^++\si_l^-})\ge  e^{-\si_l^-}/2. $$ Therefore
\begin{equation}
|\De_l(z)| \le \frac{6d_l|z|}{e^{\si_l^-}} \le  24|z|
\int\limits_{\si_l^-}^{\si_l^+} e^{-\si} \log  \psi(e^\si) \,
d\si.\label{e3.6'}
\end{equation}
Since $\psi(x)$ is slowly varying, so is $\log \psi(x)$, and
consequently
 $$ \log  \psi(x2^{k+1})\le (1+\ve) \log  \psi(x 2^k) \le
(1+\ve)^{k+1}\log  \psi(x), \quad  \ve>0, \qh x\ge x_\ve. $$
Applying the obtained estimates and the fact that every point $s$
belongs to the interiors of at most 4 of $\pl$s, from
(\ref{e3.6'}) we obtain
%\begin{gather}\nonumber
\begin{align}
&\sum_{\begin{substack} {l\in \L^- \\ \la_l>3} \end{substack}}
|\Delta_l(z)| \le 24 |z| \int\limits_{\log  2^{n+1}}^{+\infty}
e^{-\si} \log  \psi(e^{\si})\, d\si \nonumber
\\
&\qq\le  24|z| \sum_{k=0}^{+\infty} \int\limits_{\log
2^{n+1+k}}^{\log  2^{n+2+k}}e^{-\si} \log  \psi(e^{\si})\, d\si
\nonumber
\\
&\qq\le 24|z| \sum_{k=0}^{+\infty} \frac{\log
\psi(2^{n+k+2})}{2^{n+k+2}} \nonumber
\\
&\qq\le 24|z|  \log \psi (2^{n}) \sum_{k=0}^{+\infty}
\frac{(1+\ve)^{k+2}}{2^{n+2+k}} \le 96 \frac{\log
\psi(2^{n})}{1-\ve}. \label{e3.7}
%\end{gather}
\end{align}
Using the first equation from (\ref{e3.1}), it is not difficult to
obtain the following representation for~$\De_l$:
\begin{equation}
 \Delta _l(z)=\int\limits_{\pl} \Bl \log  \Bm 1-\frac {e^{\omega}}z\Bm  -\frac 12 \log  \Bm 1-\frac {e^{\om1}}z\Bm
-\frac12 \log  \Bm 1-\frac{e^{-\om2}}z\Bm\Br\, d\nl .
\label{e3.4'}
\end{equation}
We use (\ref{e3.4'}) for the estimate of $\De_l$ for $l\in \L^+$
with $\la_l>3$. For the analytic function $L_2(\w)=\log  \big(1-
\frac {e^\w}z\big)$ and $z\in \pl$, $l\in \L^+$ relationships
(\ref{e3.5}) hold with $L=L_2$. Therefore, similarly to
(\ref{e3.6}), we obtain the estimate
 $$ |\Delta_l(z)| \le \frac
{3d_l}{\inf\limits_{s\in \pl} |1- ze^{-s}|}\le \frac {12
e^{\si_l^+} d_l}{|z|}, \quad l\in \L^+. $$ Whence, using the
estimate $d_l\le \log  \psi(2^{n-1})$ ($l\in\L^+$), we deduce
%\begin{gather}
\begin{align}
\nonumber &\sum_{\begin{substack} {l\in \L^+ \\  \la_l>3}
\end{substack}} |\Delta_l(z)| \le C_3 \frac{\log
\psi(2^{n-1})}{|z|} \sum_{\begin{substack} {l\in \L^+ \\  \la_l>3}
\end{substack}} \int\limits_{\si_l^-}^{\si_l^+} e^\si\, d\si
\nonumber
\\
&\qq\le 4C_3 \frac{\log  \psi(2^{n-1})}{|z|}  \int\limits_{\log
R_1} ^{\log  2^{n-1}} e^\si\, d\si \le C_4 \log  \psi(2^{n-1}).
\label{e3.8}
%\end{gather}
\end{align}
For $l\in \L^+\cup\L^-$, it remains to estimate $\De_l$ for which
$1 \le  \la_l\le 3$. Then $d_l^2 \asymp m(\pl)$, in particular,
$d_l\le 6\pi$. In this case, following \cite{LM} and using the
second equality from (\ref{e3.5}), we obtain
%\begin{gather*}
\begin{align*}
\Delta_l(z) &= \Re \int\limits_{\pl} \biggl  ((\w-\om1)L'(\om1)-\frac12 L'(\om1)
(\om2 -\om1)
\\
&\qq\qq\qu+\int\limits_{\om1}^{\omega} L''(s)(\w-s)\, ds -\frac12
\int\limits_{\om1}^{\om2} L''(s)(\om2-s)\, ds \biggr)\, d\nl
\\
&=\Re \biggl\{ -z \int\limits_{\pl}
\int\limits_{\om1}^\w \frac {e^s(\w-s)}{(e^s-z)^2}\, ds
\, d\nl
\\
&\qq\qu+\frac z2 \int\limits_{\pl}
\int\limits_{\om1}^{\om2} \frac {e^s(\om2-s)}{(e^s-z)^2}\, ds\, d\nl\biggr\} .
%\end{gather*}
\end{align*}
Whence,
\begin{equation} \label{e3.9}
|\De _l(z)|\le 6 d_l^2 |z| \sup_{s\in \pl} \Bigl| \frac{e^s}{(e^s-z)^2} \Bigr|.
\end{equation}
If $l\in \L^+$, we have $|e^s -z|\ge |z|/2$ and from (\ref{e3.9})
we obtain
\begin{equation} \label{e3.10}
|\De_l(z)| \le 24 d_l^2 \frac{e^{\si_l^+}}{|z|} \le  \frac{C_4}{|z|} \iint\limits_{\pl}
e^\si\, d\si\,dt.
\end{equation}
Applying (\ref{e3.10}), we deduce
%^^@(23.03.04 20:31)
%\begin{gather}\nonumber
\begin{align}
&\sum_{\begin{substack} {l\in \L^+ \\ \la_l\le 3} \end{substack}}
|\De_l(z)|
%\nonumber\\
\le \frac{C_4}{|z|} \sum_{\begin{substack} {l\in \L^- \\ \la_l\le 3} \end{substack}}
\iint\limits _{\pl} e^\si\, d\si\,dt  \nonumber\\ &\qq
\le \frac{4C_4}{|z|} \int\limits_0^{2\pi}
%русское С:
%^^@\int\limits _{\log  R_1}^{\log  2^{n-2}}  e^\si\, d\si\,dt \le
%%\frac {8\pi С_4 2^{n-2}}{|z|} = O(1).
%латинское С (правильно):
\int\limits _{\log  R_1}^{\log  2^{n-2}}  e^\si\, d\si\,dt
\le
\frac {8\pi C_4 2^{n-2}}{|z|} = O(1).
\label{e3.11}
%\end{gather}
\end{align}
Now, let $l\in \L^-$, then $|e^s-z|\ge e^\si/2$, and therefore,
from (\ref{e3.9}) we obtain $$ |\De_l(z)| \le  24 d_l^2
|z|e^{-\si_l^-} \le C_5 |z| \iint\limits_{\pl} e^{-\si}\,d\si dt,
$$
\begin{equation*}%\label{e3.12}
%^^@(23.03.04 20:44)
%^^@\aligned
\sum_{l\in \L^-, \la_l\le 3} |\De_l(z)|
\le  8\pi K_5|z| \int\limits_{\log  2^{n+1}}^{\infty }
e^{-\si} \, d\si= 8\pi |z|2^{-n-1}
\le C_6.
%\endaligned
\end{equation*}
Taking into account (\ref{e3.7}), (\ref{e3.8}), (\ref{e3.11}) and
the latter estimate, we obtain
\begin{equation} \label{e3.12}
\sum_{l\in \L^+ \cup \L^-} |\De_l(z)|\le C_7 \log  \psi (2^{n+1}).
\end{equation}
Now, let $l\in \L^0$, i.~e. $Q^{(l)} \cap G_n\ne \varnothing $,
where
%^^@(23.03.04 20:44)
$G_n=\{z: 2^{n-2} \le |z| \le^^@\break 2^{n+1}\}$.
%^^@(23.03.04 20:44)
{\tolerance4000\looseness1 Among $l\in \L^0$, there are at most 8
such that $\la_l
>3$.
Indeed, for such  $l$ we have $Q^{(l)}=\{z: e^{\si_l^-} \le |z|
\le e^{\si_l^+}\}$ and $\si_l^+-\si_l^-> 6\pi$, therefore,
$e^{\si_l^+}-e^{\si_l^-}\ge 2^{n-1} (1-e^{-6\pi})$. Since the
interiors of these $Q^{(l)}$ form at most 4-fold cover, for every
$l$, excepting at most 4, we have  $e^{\si_l^-}\ge 2^{n-1}$, and,
excepting at most 8, $e^{\si_l^-}\ge e^{6\pi} 2^{n-1}> 2^{n+1}$,
i.e. the intersection with $G_n$ is empty.\par}

%^^@(23.03.04 20:44)
{\tolerance4000\looseness1 We denote this exceptional set of
indices $l$ by $\L_*^0$. We have to estimate
 $\int_{G_n}
|\De_l(z)| \, dm(z)$ for $l\in \L_*^0$. If $|\zeta|\ge 3\cdot 2^n$
($\zeta \in Q^{(l)}$, $z\in G_n$), then $\big| \log  | 1-\frac
z\zeta |\big| \le 2|z|/|\zeta|\le 3$. Otherwise, we have $|\zeta|
\le ^^@{3\cdot 2^n}$,
%\begin{gather}\nonumber
\begin{align}
&\int\limits_{D_\zeta(2^n)} \Bm \log  \Bm 1-\frac z\zeta \Bm\Bm \,
dm(z) \nonumber\\ &\qq
\le
\int\limits_{D_\zeta(2^n)} \Bm \log  \Bm \frac {z-\zeta}{2^n}
\Bm\Bm +\Bm \log  \Bm \frac \zeta{2^n}\Bm\Bm \, dm(z) \nonumber
\\
&\qq\le 2\pi \int\limits_0^{2^n} \log   \frac{2^n}{\tau} \tau \,
d\tau +\sup_{\zeta \in Q^{(l)}} \Bm \log  \Bm  {\zeta}2^{-n} \Bm
\Bm \pi 2^{2n} \nonumber\\ &\qq \le \pi^2 4^{n-1} + C_8   \log
\psi(2^n)\cdot 4^n. \label{e3.11'}
%\end{gather}
\end{align}
\par}
\noindent Thus,
\begin{align}\allowdisplaybreaks\nonumber
&\int\limits_{G_n} \sum_{l\in \L^*} \biggl| \int\limits_{Q^{(l)}}
\log  \Bigl| 1-\frac z\zeta \Bigr| \, d\mu^{(l)}(\zeta)\, \biggr|
dm(z)
\\ \allowdisplaybreaks
\nonumber &\qq\le  \sum_{l\in \L^*} \biggl(
\int\limits_{Q^{(l)}\cap \{|\zeta|\ge 3\cdot 2^n\}}+
\int\limits_{Q^{(l)}\cap  \{|\zeta|\le 3\cdot 2^n\}} \biggr)
\int\limits_{G_n} \Bm \log  \Bm 1-\frac z\zeta \Bm\Bm \, dm(z)
d\mu^{(l)}(\zeta)
\\ \allowdisplaybreaks
\nonumber
&\qq\le  \sum_{l\in \L^*} \biggl(3\cdot 2 m(G_n)
\\ \allowdisplaybreaks
\nonumber &\qq\qq\qq \aligned + \int\limits_{Q^{(l)}\cap
\{|\zeta|\le 3\cdot 2^n\}} \biggl(&\int\limits_{D_\zeta(2^n)} \Bm
\log  \Bm 1-\frac z{\zeta} \Bm \Bm \, dm(z)
\\
&+\int\limits_{G_n \setminus D_\zeta(2^n)} \Bm \log  \Bm 1-\frac
z\zeta\Bm \Bm \, dm(z)\biggl)d\mu^{(l)}(\zeta)\biggr)
\endaligned
\\
\allowdisplaybreaks \nonumber &\qq\le  K_9 4^n + K_9 \log
\psi(2^n) \cdot 4^n
\\ \allowdisplaybreaks
&\qq\qu+ \int\limits_{Q^{(l)}\cap  \{|\zeta|\le 3\cdot 2^n\}}
\int\limits_{G_n \setminus D_\zeta(2^n)} \log
\frac{|z|+|\zeta|}{|\zeta|} \, dm(z)\, d\mu^{(l)} (\zeta)
\nonumber\\ &\qq\le K_9 4^n \log  \psi(2^n). \label{e3.12'}
\end{align}
For $l\in \L^0\setminus \L^0_*$ we have $\la_l\le 3$, i.e. all the
corresponding $\pl$s are ``almost squares''; therefore, we can
apply the arguments from \cite[e.--g.]{LM}
\newcommand{\dist}{\mathop{\rm dist}}
For $D_l=\diam Q^{(l)}$, under the condition $\dist \{z,
Q^{(l)}\}>4 D_l$, using the last equality from (\ref{e3.5}),
similarly as (\ref{e3.9}) we obtain
\begin{align} \nonumber
|\De_l(z)|
&=\Bm \frac12 \Re \int\limits_{\pl} z\int\limits_{\om1}^{\omega} e^s
\frac {z+e^s}{(e^s-z)^3}(\w-s)^2\, ds\, d\nl
\\
\nonumber
&\qq\qq+\frac 12 \Re z \int\limits_{\om1}^{\om2} e^s
\frac {z+e^s}{(e^s-z)^3}(\om2-s)^2\, ds\Bm
% \le
\\
&\le  \frac {d_l^3 |z|^3}{|\zeta_l^{(1)} -z|^3} \le \frac{D_l^3}{|\zeta_l^{(1)} -z|^3}.
\label{e19}
\end{align}
Whence
\begin{align}\allowdisplaybreaks\nonumber
%\begin{align} \nonumber
&\int\limits_{G_n}
\sum_{l\in \L^0\setminus \L_*^0} |\De_l(z)| \, dm(z)
\le
 \sum_{l\in \L^0, \la_l\le 3}\int\limits_{G_n} |\De_l(z)| \, dm(z)
\nonumber\\ \allowdisplaybreaks
&\qu\le \sum_{l\in \L^0, \la_l\le 3}\biggl( \int\limits_{G_n \cap \{|z-\zeta_l^{(1)}| >3D_l\}} +
\int\limits_{ |z-\zeta_l^{(1)}| \le 3D_l}\biggr) |\De_l(z)|\, dm(z)  \nonumber
\\ \allowdisplaybreaks
&\qu\le \sum_{l\in \L^0, \la_l\le 3}\biggl( \int\limits_{G_n \cap \{|z-\zeta_l^{(1)}| >3D_l\}}
\frac {D_l^3}{|z-\zeta_l^{(1)}|^3}\, dm(z) \nonumber
\\
&\qq\qq\qq\qq\qu+\int\limits_{ |z-\zeta_l^{(1)}| \le 3D_l} |\De_l(z)|\, dm(z)\biggr).
\label{e20}
\end{align}
For the first sum from (\ref{e20}) we have
\begin{align} \nonumber
&\sum_{l\in \L^0, \la_l\le 3} D_l^3 \int\limits _{ \{|z-\zeta_l^{(1)}| > 3D_l\}\cap G_n}
\frac {1}{|z-\zeta_l^{(1)}|^3}\, dm(z)
\\
&\qq\le \sum_{l\in \L^0, \la_l\le 3} D_l^3
2\pi \int\limits_{3D_l}^{2^{n+1}} \frac 1{t^3} t\,dt  \nonumber
\\
&\qq\le \sum_{l\in \L^0, \la_l\le 3} 2D_l^2 \le \sum_{l\in \L^0, \la_l\le 3}
m(Q^{(l)}).
\label{e21}
\end{align}
It remains to estimate $\int_{ |z-\zeta_l^{(1)}| \le 3D_l}
|\De_l(z)|\, dm(z)$. From the definition of $\De(z)$ and the
equality $\log  |\zeta_l^{(1)} |+\log  |\zeta_l^{(2)} |=
\int_{Q^{(l)} }\log  |\zeta|\, d\mu^{(l)}(\zeta)$ we deduce
 $$ |\De_l(z)|=
\int\limits_{Q^{(l)}} \Bl \log  \Bm \frac{\zeta -z}{3D_l} \Bm
-\frac 12 \log  \Bm\frac{\zeta_l^{(1)} -z}{3D_l}\Bm -\frac12 \log
\Bm \frac{\zeta_l^{(2)} -z}{3D_l}\Bm \Br\, d\mu^{(l)}(\zeta). $$
Thus,
\begin{align*}\allowdisplaybreaks
&\int\limits_{|z-\ze1| \le 3D_l} |\De_l(z)|\, dm(z)
\\
&\qu\le \int\limits_{Q^{(l)}} \biggl ( \int\limits_{ |z-\ze1| \le
3D_l} \Bm\log  \Bm \frac{\zeta -z}{3D_l} \Bm\Bm +\frac 12 \Bm\log
\Bm\frac{\zeta_l^{(1)} -z}{3D_l}\Bm\Bm
\\
&\phantom{\int\limits_{|z-\ze1| \le 3D_l} |\De_l(z)|\, dm(z)\,
dm(z)}+\frac12 \Bm\log  \Bm \frac{\zeta_l^{(2)} -z}{3D_l}\Bm \Bm\,
dm(z)\biggr) d\mu^{(l)}(\zeta)
\\
&\qq= \int\limits_{Q^{(l)}} \biggl(  \int\limits_{|z-\zeta|\le
D_l} +\int\limits_{\begin{substack}{ |z-\zeta|>D_l \\
|z-\zeta_l^{(1)}|\le 3D_l} \end{substack}}\biggr) \Bm\log  \Bm
\frac{\zeta -z}{3D_l} \Bm\Bm  dm(z) \,d\mu^{(l)}(\zeta)
\\
&\qq\qu+\int\limits_{|z-\zeta_l^{(1)}|<D_l}  \log
\frac{3D_l}{|\zeta_l^{(1)} -z|} dm(z)
\\ \allowdisplaybreaks
&\qq\qu+\biggl(\int\limits_{|z-\zeta_l^{(2)}|\le 3D_l}  +\int\limits_{
\begin{substack}{
|z-\zeta_l^{(2)}|>3D_l \\ |z-\zeta_l^{(1)}|\le 3D_l}
\end{substack}} \biggr) \Bm \log  \Bm \frac{\ze2 -z}{3D_l}\Bm
\Bm\, dm(z)
\\
&\qu\le \int\limits_{Q^{(l)}} \biggl(\int\limits_0^{D_l} \log
\frac {3D_l}\tau \tau\, d\tau+ (3D_l)^2 \pi \log  3\biggr) \,
d\mu^{(l)}(\zeta)
\\
&\qq+ 2 \int\limits_0^{D_l} \log  \frac{3D_l}{\tau} \tau \,d\tau +
\pi (3D_l)^2 \log  3
\\
&\qu\le  9(3\pi \log  3+1)  D_l^2 \le C_{11} m(Q^{(l)}).
\end{align*}
From the latter inequality and (\ref{e20}) we obtain
 $$
\int\limits_{G_n} \sum_{l\in \L^0\setminus \L^0_*}
|\De_l(z)|\,dm(z) \le \sum_{l\in \L^0\setminus \L^0_*}K_{11}
m(Q^{(l)}) \le K_{12} 2^{2n}. $$
Applying the latter inequality
and (\ref{e3.12}) completes the proof of Theorem~1.
%\end{proof}
Using Chebyshev's inequality, from (\ref{e2.2}) one can easily
obtain Corollary~1.
\section{Proof of Theorem $\bf 2'$}
Suppose that $\si$ satisfies the conditions of the theorem.
Without loss of generality, one may assume that
$\psi(r)=\exp\bigl\{\int_1^r \frac{\si(t)}t \, dt\bigr\}$ is
unbounded. Obviously,  $\psi \in\Phi$. Let  $u_\psi$ be defined by
formula (\ref{ex2.1}) with $\vfi=\psi$, $\psi\in \Phi$. Suppose
that there exists an entire function $f$ and a constant $\al\in
[0,1)$ satisfying the condition
%\begin{equation} \label{e4.1}
\begin{align}
&\int\limits_{|z|<R} \big|u_\psi(z)-\log  |f(z)|-\al \log
|z|\big|\, dm(z) \nonumber
\\
\label{e4.1} &\qquad< \ve R^2 \log  \psi(R), \quad R\ge R_\ve,
%\end{equation}
\end{align}
for arbitrary $\ve > 0$. Without loss of generality, one may
assume that $f(0)\ne 0$. Separating from $u_\psi$ the term $\frac
12 \log | 1-z/r_1|$, the case $\al \in[1/2, 1)$ can be reduced to
the case $\al \in [0, 1/2)$. Therefore, we consider the latter one
in details. Define the counting Nevanlinna characteristics
 of the Riesz masses
$u_\psi$ and $f$:
%$$
\begin{align*}
N(r,u_\psi)  =\int\limits_0^r \frac{n(t,u_\psi)}t \, dt, \quad
N(r,f) =\int\limits_0^r \frac{n(t,f)}t \, dt,
\end{align*}
%$$
where $n(r,f)$ is the number of zeros of $f$ in
$\overline{D}_0(r)$. By the Jensen formula \cite[Chapter
3.9]{Hay}, we have
%$$
\begin{align*}
&\frac 1{2\pi} \int\limits_0^{2\pi} (u_\psi(re^{i\theta})- \log
|f(re^{i\theta})|-\al\log  r)\, d\theta
\\
&\qq=N(r,u_\psi)-N(r,f)-\al\log  r - \log  |f(0)|.
\end{align*}
%$$
If $R\ge \widetilde R_\ve$ taking into account slow varying of
$\log \psi(R)$ we obtain
\begin{align*}
&\int\limits_R^{2R} \big|N(r,u_\psi)-N(r,f)-\al\log  r -\log
|f(0)|\big| r\,dr
\\
&\qq\le \int\limits_R^{2R}\frac 1{2\pi} \int\limits_0^{2\pi}
\big|u_\psi (re^{i\theta}) -\log  |f(re^{i\theta})|-\al \log
r\big|\ d\theta  r\, dr
\\
&\qq\le \frac 1{2\pi} \int\limits_{|z|\le 2R} \big|u_\psi(z)-\log
|f(z)|-\al \log  |z|\big|\, dm(z)
\\
&\qq<\frac{\ve}{2\pi} (2R)^2 \log  \psi(2R)
< \frac{2\ve(1+\ve)}{\pi} R^2\log \psi(R).
\end{align*}
This implies that on $[R, 2R]$ there exists $r^*$ such that for
 $\ve\in (0, 1/2)$
\begin{equation} \label{e4.2}
|N(r^*, u_\psi)-N(r^*,f)-\al \log  r^* |\le \frac {2\ve(1+\ve)}{3}
\log  \psi(r^*)< \ve \log  \psi(r^*).
\end{equation}
We are going to derive from (\ref{e4.2}) the relationship
\begin{equation} \label{e4.3}
|n(r,u_\psi) -n(r,f)-\al |\le \frac 12, \quad r\to+\infty.
\end{equation}
Assume the contrary. If (\ref{e4.3}) fails to hold, then there
exists a sequence $(\tau_k)$, $\tau_k\to+\infty$ $(k\to+\infty)$,
such that either i) $n(\tau_k, u_\psi) -n(\tau_k, f)-\al
> 1/2$ or  ii)~$n(\tau_k, f) -n(\tau_k, u_\psi)+\al > 1/2$
$(k\to+\infty)$. Consider case i). For arbitrary $t\in [\ti
\tau_k, \tau_k]$, where $\ti \tau_k \psi(\ti \tau_k)=\tau_k$, we
have $n(t,u_\psi)- n(\tau_k, u_\psi) \ge - 1/2$, therefore
 $$ n(t,
u_\psi)- n(t, f)-\al  \ge n(\tau_k, u_\psi) -n(\tau_k, f) -\al +
n(t, u_\psi) -n (\tau_k, u_\psi)  >0. $$ Since $n(t,u_\psi)
-n(t,f)$ take the values equal to an integer multiple of $1/2$, we
have
\begin{equation}
\label{e4.4'}
n(t, u_\psi)- n(t, f)-\al  \ge 1/2-\al>0, \quad t\in [\tilde \tau_k, \tau_k].
\end{equation}
Choose $t_k \in  [\ti\tau_k, 2\ti\tau_k]$, $T_k \in [\tau_k/2,
\tau_k]$ so that (\ref{e4.2}) hold for $r^*\in \{t_k,T_k\}$. Then,
by the definition of the function $N(r, \cdot)$, applying
(\ref{e4.2}) and (\ref{e4.4'}), for $\ve\in (0, 1/4 -\al/2)$ we
obtain
\begin{align*}
\ve \log  \psi(T_k)
&\qu> |N(T_k , u_\psi) -N(T_k,f)- \al\log  T_k |
\\
&\qu\ge \int\limits_{t_k}^{T_k} \frac {n(t,u_\psi)-n(t,f)-\al }t
dt   - |N(t_k , u_\psi) -N(t_k,f)-\al \log  t_k |
\\
&\qu\ge \Bigl(\frac 12 -\al \Bigr) \log  \frac {T_k}{t_k} -\ve
\log \psi(t_k)
\\
&\qq=\Bigl(\frac 12 -\al  -\ve\Bigr)\log  \psi(t_k)>\ve\log
\psi(T_k), \quad k\to+\infty.
\end{align*}
Therefore, case i) is impossible. Similarly, in case ii) we have
$n(t,f)-n(t, u_\psi)+\al > 0$  for $t\in [\tau_k, \tau_k\psi(
\tau_k)]$, and consequently $n(t,f)-n(t, u_\psi)+\al \ge \be$,
where $\be$ is a positive constant. Choosing $t_k\in [\tau_k,
2\tau_k]$, $T_k\in [\tau_k\psi(\tau_k)/2, \tau_k\psi(\tau_k)]$
satisfying  (\ref{e4.2}) instead of $r^*$, and  $\ve \in (0,
\be/2)$, similarly as before, we obtain
\begin{align*}
&|N(T_k,f) - N(T_k , u_\psi)+ \al \log  T_k|
\\
&\qu\ge
 \int\limits_{t_k}^{T_k} \frac {n(t,f) -n(t,u_\psi)+\al }t
dt - |N(t_k,f)-  N(t_k , u_\psi) +\al \log  t_k |
\\
&\qu\ge(\be-\ve)\log  \psi (t_k).
\end{align*}
Therefore, case ii) is also impossible. Thus, (\ref{e4.3}) holds.

Let $\rho_k$ be the modulus of the $k$-th zero of $f$ (the zeros
are ordered by non-decreasing of their moduli). Since the jumps of
$n(t,f)$ take the natural values and the jumps of $n(t, u_\psi)$
the value $\frac12$, relationship (\ref{e4.3}) is possible only in
the case when, starting from some $k_0\in\N$, between every two
immediate jump points $\rho_k\le \rho_{k+1}$ of the function
$n(t,f)$, there are points $r_m, r_{m+1}$, in particular, $\rho_k
<\rho_{k+1}$. We first consider the case $\al=0$. If
$\rho_k<r_{2k-1}$, then for $r\in (\rho_k, r_{2k-1})$ we have
$n(r,f)-n(r,u_\psi)\ge k- (2k-2)/2 =1$. If $\rho_k>r_{2k} $, then
$n(r,u_\psi)- n(r,f)\ge 1$ for $r\in (\max\{r_{2k},\rho_{k-1}\} ,
\rho_k)$. Therefore, $r_{2k-1} \le \rho_k \le r_{2k}$, starting
from some $k\ge k_1$. Thus,
%^^@неправильно:
%$$ n(t,u_\psi) -n(t,f)= \cases \frac 12, & t\in[r_{2k-1}, \rho_k) \\
%-\frac 12, & t\in [\rho_k, r_{2k})
%\endcases , \quad k\ge k_1.$$
%правильно (надо \begin-, \end{cases}):
$$
n(t,u_\psi) -n(t,f)
= \begin{cases} \frac 12, & t\in[r_{2k-1}, \rho_k) \\
               -\frac 12, & t\in [\rho_k, r_{2k})
\end{cases}, \quad k\ge k_1.
$$ If $\rho_k\in[r_{2k-1}, \sqrt{r_{2k-1}r_{2k}}]$,  then choosing
$r_k^*\in[\rho_k, 2\rho_k]$, $t_k^{*}\in[r_{2k}/2, r_{2k}]$, for
which  (\ref{e4.2}) holds with $r^*\in \{r_k^*, t_k^*\}$, we
obtain
\begin{align*}
&N(t_k^{*}, f)- N(t_k^*, u_\psi)
\\
&\qu=\int\limits_{r_k^*}^{t_k^*} \frac{n(s,f)-
n(s,u_\psi)}s\, ds+ N(r_k^*,f)- N(r_k^*,u_\psi)
\\
&\qu\ge \frac 12 \log \frac{t_k^*}{r_k^*} -\ve \log  \psi(r_k^*)
\ge \Bigl(\frac 12+o(1)\Bigr) \log
\frac{t_k^*}{\sqrt{r_{2k-1}r_{2k}}} -\ve\log  \psi(r_k^*)
\\
&\qu\ge \Bl \frac14 -\ve +o(1)\Br \log  \psi(t_k^*), \quad
k\to+\infty,
\end{align*}
which contradicts to (\ref{e4.2}). In the case $\rho_k\in
[\sqrt{r_{2k-1}r_{2k}}, r_{2k}]$ we choose
 $r_k^*\in [r_{2k-1}, 2r_{2k-1}]$, $t_k^*\in[\rho_k/2, \rho_k]$
satisfying (\ref{e4.2}). Taking into account that
$n(t,u_\psi)-n(t,f)=\fract12$ for $t\in [r_k^*, t_k^*]$, we again
come to the contradiction with (\ref{e4.2}). Consequently, in the
case $\al=0$ the theorem is proved.

Now, let $\al\in (0, 1/2)$. Relationship (\ref{e4.3}) is possible
only if the expression under the modulus takes the values $-\al$
and $\frac12 -\al$. Since $r_k$ strictly increases, and the jump
of $n(r,u_\psi)$ equals $\frac12 $ for $r=r_k$ and the jump
$n(r,f)$ equals 1 for $r=\rho_k$, we see that $\rho_k=r_{2k}$,
$k\ge k_2$. Then also
%^^@неправильно:
%\begin{equation}\label{e4.5} n(t,u_\psi) -n(t,f)= \cases 0, & t\in[r_{2k}, r_{2k+1}), \\
%\frac 12, & t\in [r_{2k+1}, r_{2k+2})
%\endcases , \quad k\ge k_2.\end{equation}
%правильно (надо \begin-, \end{cases}):
\begin{equation}\label{e4.5}
n(t,u_\psi) -n(t,f)
= \begin{cases} 0, & t\in[r_{2k}, r_{2k+1}), \\
         \frac 12, & t\in [r_{2k+1}, r_{2k+2})
\end{cases}, \quad k\ge k_2.
\end{equation}
Choosing $t_k^*\in [r_{2k-1}, 2r_{2k-1}]$, $r^*_k\in[r_{2k}/2,
r_{2k}]$ satisfying (\ref{e4.2}), and taking into account
(\ref{e4.5}), as above, we come to a contradiction to
(\ref{e4.2}). Therefore, there exists no entire function $f$ with
property (\ref{e4.1}) for arbitrary $\al\in[0, 1/2)$, and
consequently, for $\al\in[0, 1)$. Theorem $ 2'$ is proved.

%\medskip
The proof of Theorem 2 literally repeats that of Theorem $2'$ for
$\al=0$, with the distinction that $\psi\in \Phi$ is given by the
condition of the theorem.
%\medskip
I would like to thank to Professor O. Skaskiv who read the paper
and made valuable suggestion as well as another participants of
the Lviv inter-university seminar in the theory of analytic
functions for valuable comments that allowed to improve the
initial version of the paper. I am also indebted to the Institute
of Mathematics at the Jagellonian University for their hospitality
during my staying in Krak\'ow in April 2002, where a part of this
paper was written.
%\end{document}

%\renewcommand{\refname}{References}

\end{document}